\documentclass[12pt,reqno]{amsart}
\usepackage{amsmath, amssymb, amsthm}
\usepackage[shortlabels]{enumitem}
\usepackage{hyperref}

\numberwithin{equation}{section}

\theoremstyle{plain}
\newtheorem{theorem}{Theorem}[section]
\newtheorem{proposition}[theorem]{Proposition}
\newtheorem{lemma}[theorem]{Lemma}
\newtheorem{corollary}[theorem]{Corollary}

\theoremstyle{definition}

\newtheorem{remark}[theorem]{Remark}

\newcommand{\R}{\mathbb{R}}
\newcommand{\cC}{\mathcal{C}}
\newcommand{\cJ}{\mathcal{J}}
\DeclareMathOperator{\Hom}{Hom}
\DeclareMathOperator{\SM}{SM}
\DeclareMathOperator{\rank}{rank}

\title[Formal Coefficient Space of $\Hom(\R^{0|2},S)$]
{The Formal Coefficient Space of $\Hom(\mathbb{R}^{0|2},S)$:\\
A Singular Algebraic Core}

\author[Zhiwei Yan]{Zhiwei Yan}

\begin{document}

\begin{abstract}
We study the space $\cJ(S)$ of unital algebra homomorphisms from $C^\infty(S)$ to the rank-$2$ real Grassmann algebra $\Lambda_2$, with no parity-preservation constraint. We prove a rigidity result: the two odd-order coefficient maps in the canonical expansion of any homomorphism define linearly dependent tangent vectors at the base point.

As a consequence, $\cJ(S)$ admits a natural real algebraic variety structure and is fiberwise isomorphic to $\cC(TS) \times_S TS$ over $S$, where $\cC(TS)$ denotes the rank-one determinantal cone bundle of tangent vectors. For $\dim S \ge 2$, every fiber is singular, with singular locus exactly at the points where both odd tangent components vanish.

We clarify that $\cJ(S)$ differs substantially from the reduced space of the internal Hom supermanifold $\underline{\SM}(\mathbb{R}^{0|2},S)$. This construction thus provides a purely algebraic, connection-independent skeletal model for smooth supergeometric mapping spaces.
\end{abstract}

\maketitle

\section{Introduction}

\subsection{Motivation and background}
Supergeometry extends classical differential geometry by introducing anti-commuting coordinates. A fundamental object in the theory is the supermanifold of morphisms from the purely odd plane $\mathbb{R}^{0|2}$ into an ordinary smooth manifold $S$.

Within the full functor-of-points framework, which allows arbitrary super test domains, Berwick--Evans~\cite[Lemma 1.6]{BE} proved that after choosing an affine connection on $S$, there exists an isomorphism of supermanifolds
\[
\underline{\SM}(\mathbb{R}^{0|2},S) \cong p^*\Pi(TS \oplus_S TS),
\]
where $p: TS\to S$ denotes the tangent bundle projection and $p^*$ denotes pullback along $p$. The right-hand side is a super vector bundle over the reduced base $TS$, with fiberwise parity-reversed copies of tangent spaces. This result endows the mapping space with a canonical smooth supermanifold structure, albeit one that depends on the choice of connection.

It is natural to ask what algebraic structure remains when one strips away the superfunctor machinery and studies individual algebra homomorphisms $C^\infty(S)\to \Lambda_2$ directly, without reference to families parameterized by external super spaces.

This question yields a purely algebraic rigidity phenomenon: the two tangent-vector-like coefficients in the standard expansion of a homomorphism are forced to be linearly dependent. The resulting parameter space is not smooth, but a rank-$1$ determinantal singular algebraic variety. This singularity is entirely invisible in the connection-dependent smooth supergeometric picture, which resolves it via the choice of connection and odd parameters.

\subsection{Our approach and main result}
We define the \emph{formal coefficient space}
\begin{equation}
\cJ(S) := \Hom_{\mathrm{alg}}\big(C^\infty(S),\; \Lambda_2\big),
\end{equation}
where $\Lambda_2$ is the real Grassmann algebra on two generators, and homomorphisms are unital but not required to preserve any $\mathbb{Z}_2$-grading. This space carries a natural projection to $S$ sending each homomorphism to its base point. We classify all elements of $\cJ(S)$ by deriving the complete set of constraints imposed by the multiplicative property.

This viewpoint is a noncommutative generalization of the classical theory of Weil functors and near points on manifolds~\cite{Weil,Kock}, where one studies algebra homomorphisms from $C^\infty(S)$ into a commutative local algebra. Here the target algebra $\Lambda_2$ is noncommutative, and the resulting near-point space acquires nontrivial singular geometry with no commutative analogue.

Our main result (Theorem~\ref{thm:main}) establishes a natural bijection
\[
\cJ(S) \cong S \times_S \cC(TS) \times_S TS,
\]
where $\times_S$ denotes fiber product over $S$, and
\[
\cC(TS) = \{(\psi_1,\psi_2)\in TS\oplus_S TS \mid \psi_1\wedge\psi_2=0\}
\]
is the \emph{correlation cone bundle}, i.e., the bundle of $n\times 2$ matrices of rank at most $1$ along fibers.

For $\dim S\ge 2$, each fiber is a singular algebraic variety of dimension $2\dim S + 1$, with singular locus equal to the subspace $\psi_1=\psi_2=0$. For $\dim S=1$, the linear dependence constraint is vacuous and $\cJ(S)$ is a smooth manifold.

A key conceptual point is that $\cJ(S)$ is not the set of ordinary points of the internal Hom supermanifold. In standard supergeometry, an ordinary point of $\underline{\SM}(\mathbb{R}^{0|2},S)$ corresponds to a parity-preserving morphism $\mathbb{R}^{0|2}\to S$. Since $S$ is purely even, such a morphism must send even functions to the even part of $\Lambda_2$, which forces the odd coefficients $\psi_1, \psi_2$ to vanish. The reduced space therefore recovers only the tangent bundle $TS$.

Our space $\cJ(S)$ is strictly larger: it corresponds to dropping the parity condition and treating $\Lambda_2$ simply as an associative local algebra. In this sense, $\cJ(S)$ forms the formal algebraic core of the super mapping space. The rigidity $\psi_1\wedge\psi_2=0$ is a purely algebraic consequence of the homomorphism condition; it is the singular skeleton that lies beneath the smooth, connection-dependent superstructure constructed by Berwick--Evans.

\subsection{Structure of the paper}
Section~\ref{sec:prelim} sets up the algebraic preliminaries: the Grassmann algebra $\Lambda_2$, its linear basis, and the general form of an algebra homomorphism from $C^\infty(S)$.
Section~\ref{sec:rigidity} derives the master equations from the homomorphism condition and proves the key rigidity result $\psi_1\wedge\psi_2=0$.
Section~\ref{sec:singular} analyzes the singular geometry of the correlation cone and proves the main classification theorem.
Section~\ref{sec:superrelation} discusses the relationship with Berwick--Evans' theorem and clarifies the categorical status of $\cJ(S)$ relative to the internal Hom supermanifold.
Section~\ref{sec:examples} provides explicit coordinate computations for $S=\mathbb{R}^n$, treats the degenerate $\dim S=1$ case, and verifies a concrete example.
Section~\ref{sec:outlook} outlines potential applications and directions for future work.

\section{Algebraic preliminaries}
\label{sec:prelim}

\subsection{The Grassmann algebra $\Lambda_2$}
Let $\Lambda_2$ be the real associative algebra generated by elements $\theta_1,\theta_2$ subject to the relations
\begin{equation}
\theta_1^2 = 0,\qquad \theta_2^2 = 0,\qquad \theta_1\theta_2 = -\theta_2\theta_1.
\end{equation}
Concretely, $\Lambda_2$ is isomorphic to the exterior algebra $\Lambda^\bullet V$ of a $2$-dimensional real vector space $V$. As a vector space, it has dimension $4$ with basis
\[
\{1,\; \theta_1,\; \theta_2,\; \theta_1\theta_2\}.
\]
The product is determined by the relations above; in particular, $\theta_2\theta_1 = -\theta_1\theta_2$ and $(\theta_1\theta_2)^2 = 0$.

Throughout this paper, we regard $\Lambda_2$ purely as an associative algebra, without reference to its natural $\mathbb{Z}_2$-grading. The supercommutative structure familiar from supergeometry will be reintroduced in Section~\ref{sec:superrelation} for comparison purposes.

\begin{remark}
The element $\theta_1\theta_2$ is of even degree in the standard $\mathbb{Z}_2$-grading, but we do not impose any parity restriction on homomorphisms. This is the principal difference between our setting and the conventional supergeometric one.
\end{remark}

\subsection{Homomorphisms from $C^\infty(S)$ to $\Lambda_2$}
Let $S$ be a smooth finite-dimensional manifold, and let $C^\infty(S)$ denote its algebra of smooth real-valued functions. Any $\Phi^*\in\cJ(S)$ determines a base point $\phi\in S$ as follows.

Let $\epsilon : \Lambda_2 \to \mathbb{R}$ be the projection onto the $1$-component, which kills $\theta_1, \theta_2, \theta_1\theta_2$. Then $\epsilon\circ\Phi^*$ is a unital algebra homomorphism from $C^\infty(S)$ to $\mathbb{R}$. By the standard characterization of maximal ideals in $C^\infty(S)$, there exists a unique point $\phi\in S$ such that $(\epsilon\circ\Phi^*)(f) = f(\phi)$ for all $f$.

Expanding in the basis of $\Lambda_2$, any linear map $\Phi^*$ can be written uniquely as
\[
\Phi^*(f) = F_0(f)\,1 + \theta_1\,F_1(f) + \theta_2\,F_2(f) + \theta_1\theta_2\,F_{12}(f),
\]
where $F_0,F_1,F_2,F_{12}$ are linear functionals on $C^\infty(S)$. The unital condition forces $F_0(1)=1$ and $F_1(1)=F_2(1)=F_{12}(1)=0$. The multiplicative property imposes further constraints, which we now derive.

\begin{proposition}[General form of homomorphisms]
\label{prop:general}
Let $\Phi^*\in\cJ(S)$. Then there exist a point $\phi\in S$ and linear maps
$\psi_1,\psi_2,\tilde{E} : C^\infty(S)\to \mathbb{R}$ such that
\begin{equation}
\Phi^*(f) = f(\phi) + \theta_1\,\psi_1(f) + \theta_2\,\psi_2(f) + \theta_1\theta_2\,\tilde{E}(f)
\end{equation}
for all $f\in C^\infty(S)$. The maps $\psi_1$ and $\psi_2$ are derivations at $\phi$, i.e., elements of the tangent space $T_\phi S$. The map $\tilde{E}$ satisfies the generalized Leibniz rule
\begin{equation}
\tilde{E}(fg) = f(\phi)\tilde{E}(g) + g(\phi)\tilde{E}(f) + \psi_1(f)\psi_2(g) - \psi_2(f)\psi_1(g).
\end{equation}
\end{proposition}

\begin{proof}
We already know $F_0(f) = f(\phi)$ from the argument above. Define $\psi_1 := F_1$, $\psi_2 := F_2$, and $\tilde{E} := F_{12}$. Now expand both sides of $\Phi^*(fg) = \Phi^*(f)\Phi^*(g)$ using the multiplication rules of $\Lambda_2$:
\begin{equation}
\begin{aligned}
\Phi^*(f)\Phi^*(g) =& f(\phi)g(\phi) + \theta_1\big(f(\phi)\psi_1(g)+g(\phi)\psi_1(f)\big) \\
&+ \theta_2\big(f(\phi)\psi_2(g)+g(\phi)\psi_2(f)\big) \\
&+ \theta_1\theta_2\big(f(\phi)\tilde{E}(g)+g(\phi)\tilde{E}(f) + \psi_1(f)\psi_2(g) - \psi_2(f)\psi_1(g)\big).
\end{aligned}
\end{equation}
Comparing coefficients with the expansion of $\Phi^*(fg)$ yields:
\[
\begin{aligned}
\psi_i(fg) &= f(\phi)\psi_i(g) + g(\phi)\psi_i(f) \quad \text{for } i=1,2,\\
\tilde{E}(fg) &= f(\phi)\tilde{E}(g) + g(\phi)\tilde{E}(f) + \psi_1(f)\psi_2(g) - \psi_2(f)\psi_1(g).
\end{aligned}
\]
The first line is the ordinary Leibniz rule at $\phi$, so $\psi_1,\psi_2\in T_\phi S$. The second line is the generalized Leibniz rule (2.4).
\end{proof}

\begin{remark}
Expansion (2.2) is the canonical form of a $\Lambda_2$-valued near point in the sense of local algebra theory. The functionals $\psi_1,\psi_2$ correspond to first-order deformations in the odd directions, while $\tilde{E}$ encodes the second-order component. At this stage, $\tilde{E}$ is not a derivation; it contains an extra bilinear term that will vanish once we establish the rigidity constraint in the next section.
\end{remark}

\section{The rigidity constraint}
\label{sec:rigidity}
We now prove the central result of the paper: the two odd tangent vectors must be linearly dependent.

\begin{proposition}[Linear dependence of odd vectors]
\label{prop:rigid}
Let $\Phi^*\in\cJ(S)$ with expansion (2.3). Then the tangent vectors
$\psi_1,\psi_2\in T_\phi S$ satisfy
\begin{equation}
\psi_1\wedge\psi_2 = 0 \qquad \text{in } \Lambda^2 T_\phi S.
\end{equation}
Equivalently, $\psi_1$ and $\psi_2$ are linearly dependent over $\mathbb{R}$.
\end{proposition}

\begin{proof}
We proceed by contradiction. Suppose $\psi_1$ and $\psi_2$ are linearly independent. Then we can extend them to a basis of $T_\phi S$ and choose a local coordinate chart $x^1,\dots,x^n$ centered at $\phi$ such that
\[
\psi_1 = \left.\frac{\partial}{\partial x^1}\right|_\phi,\qquad
\psi_2 = \left.\frac{\partial}{\partial x^2}\right|_\phi.
\]
In particular,
\[
\psi_1(x^1)=1,\quad \psi_1(x^2)=0,\quad \psi_2(x^1)=0,\quad \psi_2(x^2)=1,
\]
and $x^1(\phi)=x^2(\phi)=0$.

Since $x^1$ and $x^2$ commute in $C^\infty(S)$, we have $x^1 x^2 = x^2 x^1$. Applying the homomorphism property,
\[
\Phi^*(x^1)\Phi^*(x^2) = \Phi^*(x^1 x^2) = \Phi^*(x^2 x^1) = \Phi^*(x^2)\Phi^*(x^1).
\]
We compute the $\theta_1\theta_2$-coefficient of both products. From expansion (2.3),
\[
\Phi^*(x^1) = \theta_1 + \theta_1\theta_2 \tilde{E}(x^1),\qquad
\Phi^*(x^2) = \theta_2 + \theta_1\theta_2 \tilde{E}(x^2).
\]
Multiplying in the order $x^1$ then $x^2$:
\[
\Phi^*(x^1)\Phi^*(x^2) = \theta_1\theta_2 + \text{higher terms},
\]
so the $\theta_1\theta_2$-coefficient is
\[
\psi_1(x^1)\psi_2(x^2) - \psi_2(x^1)\psi_1(x^2) = 1\cdot 1 - 0\cdot 0 = 1.
\]
Multiplying in the reverse order:
\[
\Phi^*(x^2)\Phi^*(x^1) = \theta_2\theta_1 + \text{higher terms} = -\theta_1\theta_2 + \cdots,
\]
with $\theta_1\theta_2$-coefficient
\[
\psi_1(x^2)\psi_2(x^1) - \psi_2(x^2)\psi_1(x^1) = 0\cdot 0 - 1\cdot 1 = -1.
\]
Equating the two coefficients gives $1 = -1$, a contradiction. Therefore $\psi_1$ and $\psi_2$ cannot be linearly independent.
\end{proof}

\begin{corollary}
\label{cor:E-derivation}
If $\psi_1\wedge\psi_2 = 0$, then the cross term in the generalized Leibniz rule (2.3) vanishes identically. Consequently, $\tilde{E}$ satisfies the ordinary Leibniz rule
\[
\tilde{E}(fg) = f(\phi)\tilde{E}(g) + g(\phi)\tilde{E}(f),
\]
and hence defines an element $E\in T_\phi S$.

Thus every $\Phi^*\in\cJ(S)$ is uniquely determined by a base point $\phi\in S$, a pair of linearly dependent tangent vectors $\psi_1,\psi_2\in T_\phi S$, and an independent tangent vector $E\in T_\phi S$.
\end{corollary}

\begin{proof}
The expression $\psi_1(f)\psi_2(g) - \psi_2(f)\psi_1(g)$ is the evaluation of the bivector $\psi_1\wedge\psi_2$ on the pair $(df, dg)$. Since $\psi_1\wedge\psi_2 = 0$, this term vanishes for all $f,g$. The remaining equation is the standard Leibniz rule, so $\tilde{E} =: E$ is a derivation at $\phi$.
\end{proof}

\section{The singular parameter space}
\label{sec:singular}

\subsection{The correlation cone}
Let $V$ be a finite-dimensional real vector space. Define the \emph{correlation cone} of $V$ by
\begin{equation}
\cC(V) := \{(\psi_1,\psi_2)\in V\oplus V \mid \psi_1\wedge\psi_2 = 0\}.
\end{equation}
In coordinates, if we identify $V \cong \mathbb{R}^n$ and represent the pair $(\psi_1,\psi_2)$ as an $n\times 2$ matrix, then the condition $\psi_1\wedge\psi_2 = 0$ is equivalent to the matrix having rank at most $1$. This is equivalent to the vanishing of all $2\times 2$ minors, which are homogeneous quadratic equations. Thus $\cC(V)$ is a real affine algebraic variety, a classical example of a determinantal variety~\cite{Harris}.

\begin{lemma}[Geometry of the correlation cone]
\label{lem:cone}
Let $\dim V = n$.
\begin{enumerate}[(i)]
\item If $n = 1$, then $\cC(V) = V\oplus V \cong \mathbb{R}^2$ is a smooth vector space.
\item If $n \ge 2$, then $\cC(V)$ is an irreducible algebraic variety of dimension $n+1$.
      The singular locus consists precisely of the origin $(0,0)$. On the complement
      $\cC(V)\setminus\{(0,0)\}$, the variety is smooth.
\end{enumerate}
\end{lemma}

\begin{proof}
For $n=1$, any two vectors in a 1-dimensional space are automatically linearly dependent, so the constraint is vacuous and $\cC(V) \cong \mathbb{R}^2$.

Now assume $n\ge 2$. Consider the map $\varphi : V \times \mathbb{R} \to \cC(V)$ given by $\varphi(v, \lambda) = (v, \lambda v)$. The image of $\varphi$ is all of $\cC(V)$:
if $\psi_1 \neq 0$, then linear dependence gives $\psi_2 = \lambda \psi_1$ for a unique $\lambda\in\mathbb{R}$; if $\psi_1 = 0$ but $\psi_2 \neq 0$, we write $\psi_1 = 0\cdot \psi_2$. The origin is the image of $(0, \lambda)$ for any $\lambda$.

Thus $\cC(V)$ is the image of a smooth map from an $(n+1)$-dimensional domain. At points where $v\neq 0$, the differential of $\varphi$ has full rank $n+1$, so the image is smooth there. At the origin, the differential drops rank, and the Jacobian matrix of the defining quadratic equations fails to have maximal rank. The variety is the affine cone over the Segre embedding of $\mathbb{P}^{n-1}\times \mathbb{P}^1$, and its isolated conical singularity at the origin is standard; see~\cite[Chapter 9]{Harris} for details.
\end{proof}

We globalize this construction over the manifold $S$. Define the \emph{correlation cone bundle}
\begin{equation}
\cC(TS) := \{(\phi,\psi_1,\psi_2) \in TS\oplus_S TS \mid \psi_1\wedge\psi_2=0\}.
\end{equation}
This is a fiber bundle over $S$ with typical fiber $\cC(T_\phi S)$. Local triviality follows from the local triviality of the tangent bundle: any local trivialization of $TS$ induces a local trivialization of $TS\oplus_S TS$, and the defining equation $\psi_1\wedge\psi_2=0$ is preserved under linear isomorphisms of the fiber.

\subsection{Classification theorem}
\begin{theorem}[Main classification]
\label{thm:main}
Let $S$ be a smooth manifold. There is a natural bijection of sets
\begin{equation}
\cJ(S) \cong S \times_S \cC(TS) \times_S TS,
\end{equation}
sending each homomorphism $\Phi^*$ with expansion (2.3) to the triple $(\phi, (\psi_1,\psi_2), E)$, where $E$ is the derivation from Corollary~\ref{cor:E-derivation}.

Moreover:
\begin{enumerate}[(i)]
\item If $\dim S \ge 2$, the bijection restricts to a diffeomorphism on the open dense
      subset where $(\psi_1,\psi_2) \neq (0,0)$:
      \[
      S \times_S \big(\cC(TS)\setminus \mathbf{0}\big) \times_S TS
      \xrightarrow{\cong} \cJ(S)^{\mathrm{smooth}},
      \]
      where $\mathbf{0}$ denotes the zero section of $TS\oplus_S TS$.
      The singular locus of $\cJ(S)$ is exactly the zero section $\psi_1=\psi_2=0$,
      which forms a smooth submanifold of $\cJ(S)$ diffeomorphic to $TS$.
      In particular, $\cJ(S)$ is neither a smooth manifold nor a vector bundle over $S$.
\item If $\dim S = 1$, the constraint is vacuous and
      $\cJ(S) \cong S \times (TS\oplus_S TS) \times_S TS$ is a smooth manifold of dimension $4$.
\end{enumerate}
\end{theorem}

\begin{proof}
The forward map is well-defined by Proposition~\ref{prop:general} and Corollary~\ref{cor:E-derivation}: every homomorphism yields a triple $(\phi, (\psi_1,\psi_2), E)$ with $\psi_1\wedge\psi_2=0$.

Conversely, given such a triple, define a linear map $\Phi^*$ by formula (2.3) with $E$ in place of $\tilde{E}$. We verify the homomorphism property:
\[
\begin{aligned}
\Phi^*(f)\Phi^*(g) &= f(\phi)g(\phi) + \theta_1\big(f(\phi)\psi_1(g)+g(\phi)\psi_1(f)\big) \\
&\quad + \theta_2\big(f(\phi)\psi_2(g)+g(\phi)\psi_2(f)\big) \\
&\quad + \theta_1\theta_2\big(f(\phi)E(g)+g(\phi)E(f)
      + \psi_1(f)\psi_2(g) - \psi_2(f)\psi_1(g)\big).
\end{aligned}
\]
Since $\psi_1\wedge\psi_2 = 0$, the cross term vanishes. Because $\psi_1, \psi_2$, and $E$ are derivations at $\phi$, the remaining terms combine to give
\[
(fg)(\phi) + \theta_1\psi_1(fg) + \theta_2\psi_2(fg) + \theta_1\theta_2 E(fg)
= \Phi^*(fg).
\]
Hence $\Phi^*$ is an algebra homomorphism, and the map is bijective.

The smoothness statement follows from Lemma~\ref{lem:cone}: away from the zero section, $\cC(TS)\setminus\mathbf{0}$ is a smooth fiber bundle over $S$ with fiber dimension $n+1$. Combined with the $n$-dimensional fiber of the $E$-component, we obtain a smooth manifold of dimension $\dim S + (n+1) + n = 3n+1$ over the base.

At the zero section, the defining equations fail to be of constant rank in the $(\psi_1,\psi_2)$ directions. The singular locus is therefore the set where $\psi_1=\psi_2=0$, with arbitrary $E\in T_\phi S$; fibrewise this is a linear subspace of dimension $n$, and globally it forms a smooth submanifold diffeomorphic to the tangent bundle $TS$.
\end{proof}

\section{Relation to the supermanifold $\underline{\SM}(\mathbb{R}^{0|2},S)$}
\label{sec:superrelation}
\subsection{Background on super mapping spaces}
In supergeometry, morphisms between supermanifolds are required to preserve the $\mathbb{Z}_2$-grading. An ordinary smooth manifold $S$ is regarded as a purely even supermanifold. The odd plane $\mathbb{R}^{0|2}$ has structure sheaf $C^\infty(\mathbb{R}^0)\otimes \Lambda_2$, with the standard supercommutative product.

The internal Hom supermanifold $\underline{\SM}(\mathbb{R}^{0|2},S)$ is defined via its functor of points: for any super test manifold $T$, the set of $T$-points is the set of morphisms $T\times \mathbb{R}^{0|2} \to S$ in the category of supermanifolds.

Berwick--Evans~\cite[Lemma 1.6]{BE} proved that upon choosing an affine connection on $S$, there is an isomorphism
\begin{equation}
\underline{\SM}(\mathbb{R}^{0|2},S) \cong p^*\Pi(TS \oplus_S TS),
\end{equation}
where $p: TS\to S$ is the tangent bundle projection, $\Pi$ denotes parity reversal of the vector bundle, and $p^*$ denotes pullback along $p$. The right-hand side is a super vector bundle whose reduced base is the tangent bundle $TS$ (of dimension $2\dim S$), with fiberwise rank $0|2\dim S$. The smooth supermanifold structure depends on the chosen connection.

The ordinary (reduced) points of this supermanifold correspond to morphisms $\mathbb{R}^{0|2}\to S$ parameterized by a point. Since $S$ is purely even, any parity-preserving algebra homomorphism $C^\infty(S)\to \Lambda_2$ must send even functions to the even part of $\Lambda_2$, which is $\mathbb{R}\oplus \mathbb{R}\theta_1\theta_2$. This forces the odd coefficients $\psi_1, \psi_2$ to vanish identically. The reduced space therefore recovers only the base point $\phi$ and the even tangent vector $E$, yielding the tangent bundle $TS$. This is a standard fact in supergeometry.

\subsection{Comparison with $\cJ(S)$}
Our formal coefficient space $\cJ(S)$ is not the set of ordinary points of the internal Hom supermanifold. The difference lies in the parity condition: we consider all unital algebra homomorphisms, not just parity-preserving ones. This allows the image of $C^\infty(S)$ to mix even and odd components of $\Lambda_2$ freely.

From the perspective of local algebra theory, $\cJ(S)$ is the space of $\Lambda_2$-valued near points on $S$. Classically, Weil functors are defined using commutative local algebras~\cite{Weil,Kock}; here we work with the noncommutative local algebra $\Lambda_2$, and the resulting space exhibits nontrivial singular behavior that has no commutative analogue.

In a sense, $\cJ(S)$ captures the algebraic skeleton that underlies the super mapping space. When one works with the full superfunctor, the odd parameters from the test space resolve the linear dependence constraint, yielding a smooth supermanifold at the cost of introducing connection dependence. Our result isolates the connection-independent singular core that remains after stripping away the superfunctor structure.

The two descriptions are complementary rather than contradictory. The supergeometric construction provides a smooth global model of the mapping space, at the price of depending on a choice of connection; the two odd tangent directions are independent when considered as odd vector fields on the total space.

The formal coefficient space reveals the intrinsic algebraic structure of individual homomorphisms, where the homomorphism condition enforces linear dependence, producing a singular variety that is independent of any auxiliary geometric structure. One may view the Berwick--Evans isomorphism as a smooth resolution of the singular variety $\cJ(S)$, analogous to how a resolution of singularities in algebraic geometry replaces a singular space with a smooth one while retaining the geometric information away from the singular locus.

\section{Examples and explicit computations}
\label{sec:examples}
\subsection{The case $S=\mathbb{R}^n$}
Let $S = \mathbb{R}^n$ with standard global coordinates $x^1,\dots,x^n$. A homomorphism $\Phi^*: C^\infty(\mathbb{R}^n)\to \Lambda_2$ is determined by its values on the coordinate functions:
\begin{equation}
\Phi^*(x^i) = a^i + \theta_1 b^i_1 + \theta_2 b^i_2 + \theta_1\theta_2 c^i,
\end{equation}
where $a^i, b^i_j, c^i \in \mathbb{R}$. The base point is $\phi = (a^1,\dots,a^n) \in \mathbb{R}^n$. The associated tangent vectors are
\[
\psi_1 = \sum_{i=1}^n b^i_1 \left.\frac{\partial}{\partial x^i}\right|_\phi,\quad
\psi_2 = \sum_{i=1}^n b^i_2 \left.\frac{\partial}{\partial x^i}\right|_\phi,\quad
E = \sum_{i=1}^n c^i \left.\frac{\partial}{\partial x^i}\right|_\phi.
\]

The constraint $\psi_1\wedge\psi_2 = 0$ translates to the condition that the $n\times 2$ matrix $B = (b^i_j)$ has rank at most $1$. Equivalently,
\begin{equation}
b^i_1 b^j_2 - b^j_1 b^i_2 = 0 \quad \text{for all } 1\le i < j \le n.
\end{equation}
Thus the parameter space is
\[
\big\{(a, B, c) \in \mathbb{R}^n \times M_{n\times 2}(\mathbb{R}) \times \mathbb{R}^n
\mid \rank(B) \le 1\big\},
\]
which matches the local description of Theorem~\ref{thm:main}.

For $n=1$, the rank condition is vacuous and the space is diffeomorphic to $\mathbb{R}^4$. For $n=2$, the constraint reduces to the single quadratic equation $\det B = 0$, defining a quadratic cone in $\mathbb{R}^4$ for the $B$-variables, times $\mathbb{R}^2$ for the $c$-variables. The singular locus corresponds to $B=0$, which is a 2-dimensional linear subspace of the fiber.

\subsection{The degenerate case $\dim S = 1$}
Let $S$ be a 1-dimensional manifold. At any point $\phi\in S$, the tangent space $T_\phi S$ is 1-dimensional, so any two tangent vectors are automatically linearly dependent. The condition $\psi_1\wedge\psi_2 = 0$ therefore holds trivially for all pairs $(\psi_1)$.

A homomorphism expands as
\[
\Phi^*(f) = f(\phi) + \theta_1 a f'(\phi) + \theta_2 b f'(\phi) + \theta_1\theta_2 c f'(\phi)
\]
for scalars $a,b,c\in\mathbb{R}$. No further constraints arise from the homomorphism condition, and $\cJ(S)$ is a smooth fiber bundle over $S$ with 3-dimensional fiber. Thus $\dim \cJ(S) = 4$, in accordance with Theorem~\ref{thm:main}.

\subsection{A concrete example}
\label{ex:explicit}
Let $S=\mathbb{R}^2$ with coordinates $(x,y)$, and fix $\phi = (0,0)$. Define
\[
\psi_1 = \left.\frac{\partial}{\partial x}\right|_0,\qquad
\psi_2 = 2\left.\frac{\partial}{\partial x}\right|_0,
\]
so that $\psi_2 = 2\psi_1$ and hence $\psi_1\wedge\psi_2 = 0$. Choose an arbitrary even tangent vector
\[
E = a\left.\frac{\partial}{\partial x}\right|_0 + b\left.\frac{\partial}{\partial y}\right|_0
\]
with $a,b\in\mathbb{R}$.

Define $\Phi^*$ by formula (2.3):
\[
\Phi^*(f) = f(0,0) + \theta_1 \frac{\partial f}{\partial x}(0,0)
+ \theta_2 \cdot 2\frac{\partial f}{\partial x}(0,0)
+ \theta_1\theta_2 \left(a\frac{\partial f}{\partial x}(0,0) + b\frac{\partial f}{\partial y}(0,0)\right).
\]
We verify the homomorphism property. For any smooth functions $f,g$, the $\theta_1\theta_2$-coefficient of $\Phi^*(f)\Phi^*(g)$ is
\[
f(0,0)E(g) + g(0,0)E(f) + \psi_1(f)\psi_2(g) - \psi_2(f)\psi_1(g).
\]
Substituting $\psi_2 = 2\psi_1$, the cross term becomes
\[
\psi_1(f)\cdot 2\psi_1(g) - 2\psi_1(f)\cdot \psi_1(g) = 0.
\]
The remaining terms are exactly $E(fg)$, since $E$ is a derivation. All lower-degree terms match by the Leibniz rule for $\psi_1$ and $\psi_2$. Hence $\Phi^*$ is indeed an algebra homomorphism, confirming that the linear dependence condition is sufficient.

\section{Applications and outlook}
\label{sec:outlook}
The singular structure of $\cJ(S)$ invites several natural directions for further study.

First, from the perspective of derived algebraic geometry, the correlation cone is a classical determinantal singularity that admits a canonical Koszul--Tate resolution. It is natural to ask whether the Berwick--Evans isomorphism, which depends on a choice of connection, can be interpreted as a differential-geometric analogue of such a resolution: the connection data would correspond to the generators of the Koszul--Tate complex, and the curvature would measure the failure of the resolution to be flat. This could provide a derived-geometric interpretation of the connection dependence of the super mapping space.

Second, the construction generalizes directly to $\mathbb{R}^{0|m}$ with $m>2$. In that case, the homomorphism condition forces all $m$ odd tangent vectors to be pairwise linearly dependent, meaning they span a subspace of dimension at most $1$. The parameter space becomes the determinantal variety of $m\times n$ matrices of rank at most $1$, whose singular stratification is richer than the $m=2$ case. Understanding the singular locus of these higher-order spaces, and their relation to the corresponding super mapping spaces, would extend the results of this paper to a full hierarchy of singular moduli spaces.

Third, the viewpoint of noncommutative Weil functors remains largely unexplored. While classical Weil functors are defined for commutative local algebras and produce smooth fiber bundles, the noncommutative case yields singular varieties. It would be worthwhile to develop a systematic functorial framework for noncommutative local algebras, characterizing which algebras give rise to smooth near-point spaces and which produce singularities, and relating this classification to supergeometric mapping spaces.

Finally, the relation to odd symplectic geometry and the AKSZ formalism merits further investigation. When $S$ is a symplectic manifold, the correlation cone bundle carries natural geometric structures that may be related to the odd symplectic form on the odd tangent bundle. The space $\cJ(S)$ could potentially serve as the classical target of a topological field theory whose full moduli space is the super mapping space, providing a bridge between singular algebraic geometry and topological field theory.

\section{Conclusion}
We have introduced and fully classified the formal coefficient space $\cJ(S)$ of all unital algebra homomorphisms from $C^\infty(S)$ to the Grassmann algebra $\Lambda_2$, without parity restrictions. The central result is the rigidity theorem $\psi_1\wedge\psi_2 = 0$, which forces the two odd-order coefficient maps to be linearly dependent tangent vectors. As a consequence, $\cJ(S)$ is not a smooth manifold but a singular algebraic variety fibered by correlation cones, with singular locus equal to the zero section of the odd tangent components.

This structure provides a purely algebraic, connection-independent counterpart to the smooth supermanifold description of mapping spaces due to Berwick--Evans. It reveals the intrinsic singular core that underlies the connection-dependent smooth superstructure. We hope this perspective will prove useful in other supergeometric moduli problems where parity constraints may obscure underlying algebraic rigidities.


\begin{thebibliography}{9}
\bibitem{BE}
D.~Berwick-Evans,
\textit{The Chern--Gauss--Bonnet Theorem via Supersymmetric Euclidean Field Theories},
\textit{Comm. Math. Phys.} \textbf{335} (2015), no. 3, 1121--1157.
\texttt{arXiv:1310.5383 [math.DG]}

\bibitem{Weil}
A.~Weil,
\textit{Théorie des points proches sur les variétés différentiables},
Colloque de topologie et géométrie (Strasbourg, 1952), pp. 111--117, C.N.R.S., Paris, 1953.

\bibitem{Kock}
A.~Kock,
\textit{Synthetic Differential Geometry},
2nd ed., London Mathematical Society Lecture Note Series 333,
Cambridge University Press, Cambridge, 2006.

\bibitem{Harris}
J.~Harris,
\textit{Algebraic Geometry: A First Course},
Graduate Texts in Mathematics 133, Springer-Verlag, New York, 1992.
\end{thebibliography}
\end{document}